\documentclass[11pt]{article}
\date{August, 7th, 2006}
\usepackage{amsmath}
\usepackage{amssymb}
\usepackage{amsthm}
\newcommand{\betti}{b}
\newcommand{\beq}{\begin{equation}}
\newcommand{\eeq}{\end{equation}}
\newcommand{\beqno}{\begin{equation*}}
\newcommand{\eeqno}{\end{equation*}}
\providecommand{\bysame}{\leavevmode\hbox to3em{\hrulefill}\thinspace}
\newcommand{\comp}{\mathbb{C}}

\newcommand{\Hc}{{\cal H}_c}

\newcommand{\ind}{{\rm ind}}

\newcommand{\n}{\mathbb{N}}
\newcommand{\nulli}{{\rm null}}
\newcommand{\q}{\mathbb{Q}}
\newcommand{\R}{\mathbb{R}}

\newcommand{\spec}{{\rm Spec}}

\newcommand{\symp}{{\rm Sp}}

\newenvironment{theorem}{\\[1ex] {\bf Theorem.}\em }{\\[1ex] \rm}
\newtheorem{pro}{Proposition}

\newtheorem{lem}{Lemma}
\newtheorem{rem}{Remark}

\newtheorem{dfn}{Definition}
\author{Hans-Bert Rademacher}
\title{The second closed geodesic on Finsler spheres
of dimension $n > 2$}
\begin{document}
\baselineskip 14pt
\maketitle
\begin{abstract}
We show the existence of at least
two geometrically distinct closed geodesics on an
$n$-dimensional sphere with a bumpy and
non-reversible
Finsler metric for $n > 2.$ 
\end{abstract}
\section{Introduction}
\label{sec:intro}
Estimates for the number of closed extremals for an
one-dimensional variational problem on $n-$spheres
have been studied extensively in the calculus
of variations. From the geometric point of view
one is looking 
for lower bounds for the number of geometrically 
distinct closed geodesics on $n-$dimensional spheres $S^n$
carrying a Finsler metric, cf. \cite{An}.
In \cite{LF} Lusternik and Fet prove the existence of 
one closed geodesic for a Finsler metric on a compact
and simply-connected
manifold based on ideas of Birkhoff. For further existence
results
a fundamental difference between Riemannian metrics
resp. reversible Finsler
metrics, and non-reversible Finsler metrics occurs.
On one hand
there exist non reversible Finsler metrics on 
$S^{2n}$ resp. $S^{2n-1}$ with only $2 n$ geometrically
distinct closed geodesics. These metrics occur in the
work of Katok, their geometry is studied by Ziller
in \cite{Zi}.
On the other hand there is no example of a reversible 
Finsler metric resp. a Riemannian metric on $S^n$
with only finitely many closed geodesics.
In fact on the two-sphere there are infinitely many geometrically
distinct closed geodesics for every Riemannian
metric due to the work of Bangert \cite{Ba2}
and Franks \cite{Fr}, as well as Hingston \cite{Hi2}.\\
In this paper we consider the following 
non-degeneracy assumption:
We call a Finsler metric {\em bumpy} if all closed
geodesics are non-degenerate, i.e. if there are no non-trivial
periodic Jacobi fields
orthogonal to the closed geodesic.
For a bumpy metric the energy functional on the space of
free loops is a Morse function with only non-degenerate
critical $S^1$-orbits.
This is a generic assumption
and the above mentioned Katok examples with only finitely
many closed geodesics are bumpy.
\\
The author showed in \cite[ch.4]{Ra1} that a 
bumpy Finsler metric
on $S^2$ has at least two geometrically distinct closed
geodesics. Instead of the non-degeneracy assumption one
can prove the existence of two closed geodesics on
$S^2$ resp. several geodesics on $S^n$ 
with length estimates also
for an open set of Finsler metrics defined 
for example in terms of
the flag curvature and the reversibility, cf. \cite{Ra2006}.
In their recent work \cite{BL}
Bangert and Long prove that for {\em every}
non-reversible Finsler metric on $S^2$ there are two closed
geodesics. From results by Hofer, Wysocki and Zehnder
in \cite{HW} one obtains the following statement: A
bumpy metric on $S^2$ 
for which the stable and unstable manifolds 
of all closed geodesics intersect transversally,
carries either
two or infinitely many geometrically distinct 
closed geodesics.
\\
For the $n$-sphere of dimension $n\ge 3$ there are 
only results for generic metrics
or metrics with curvature restrictions. Fet  proves
in \cite{Fe} that every bumpy and reversible Finsler metric 
on a compact manifold
carries at least two
geometrically distinct closed geodesics. 
For results on the existence of infinitely many closed
geodesics for generic metrics on spheres resp.
compact and simply-connected manifolds we refer
to \cite{Hi1},\cite{Ra3} and \cite[p.141]{Zi}. For further details
and references
we refer to the surveys
\cite{Ba1},\cite{Ta} and \cite{Lo2}.
\\
In this paper we prove for bumpy metrics
the following extension of 
the result by Bangert and Long to all dimensions:
\begin{theorem}
\label{thm:main}
Let a compact and simply-connected manifold
of the rational homotopy type of an
$n$-dimensional sphere $S^n, n\ge 3$ carry
a bumpy and non-reversible Finsler metric. Then
there are at least {\em two} 
geometrically distinct closed geodesics.
\end{theorem}
The main ingredients of the proof are the relation between
the average indices of closed geodesics for metrics with only
finitely many closed geodesics shown in \cite[Thm.3.1]{Ra1} 
and a detailed
analysis of the sequence of Morse indices $\ind(c^m)$  
of the
coverings $c^m$ of a prime closed geodesic $c$ using 
a formula due to Bott \cite{Bo} as well as
a careful discussion of the Morse inequalities. For
$n\ge 4$ we use in addition the following indirect argument:
If there is only one geometric closed geodesic
represented by the prime closed geodesic $c$ then
it is crucial that we are able to show in 
Proposition~\ref{pro:morse} that the
the sequence 
$\ind (c^m),m\ge 1$ of Morse indices is
monotone increasing. Using the Morse inequalities we
conclude that the difference 
$\ind\left(c^{m+1}\right)-\ind\left(c^m\right)\le 2$
for all $m\ge 1.$ But since $\ind(c)=n-1$ the
{\em common index jump theorem} due to Long and Zhu,
cf. \cite[Thm.4.3]{LZ} implies that there are 
integers
$k$ with 
$\ind\left(c^{2k+1}\right)-
\ind\left(c^{2k-1}\right)=2\,\ind(c) = 2n-2\ge6\,.$
\section{Critical Point Theory}
Closed geodesics 
on a compact manifold 
with a Finsler metric $F$ 
can be characterized as the critical points of the
energy functional
$$E: \Lambda M \rightarrow \R\,;\,E(\gamma)=
\frac{1}{2}\int_0^1 F^2\left(\gamma'(t)\right)\,dt \,.$$
Here $\Lambda M$ is the {\em free loop space} consisting
of closed $H^1-$curves 
$\gamma: S^1:=[0,1]/\{0,1\}\rightarrow M$ on 
the manifold $M.$ For the case of a Riemannian
manifold cf.
\cite{Kl}, for the case of a Finsler metric
cf. \cite[ch.5]{Ra2}, \cite[ch.2]{Ra4}.
References for facts from Finsler geometry are
\cite{BCS},\cite{CS}.
On $\Lambda M$ there is an $S^1$-action 
$ (u,\gamma) \in S^1 \times \Lambda M \mapsto u.\gamma \in \Lambda M;
u.\gamma(t)=\gamma(t+u), t \in S^1$
leaving the
energy functional invariant.
In addition there is the
mapping 
$^m: \gamma \in \Lambda M \mapsto \gamma^m \in \Lambda M\,;
\,\gamma^m(t)=\gamma(mt)\,;\,t \in S^1$ and 
$E(\gamma^m)=m^2 E(\gamma).$ 
Here $\gamma^m$ is the $m$-fold cover of $\gamma.$
A closed geodesic $c$ is called
{\em prime} if there is no closed geodesic $c_1$ and no
integer $m>1$ with $c=c_1^m.$ 
\\
We call a Finsler metric $F$ {\em reversible,}
if
for all tangent vectors $X$ we have: $F(-X)=F(X).$
Otherwise we call the metric {\em non-reversible.}
We call two closed geodesics $c_1,c_2:S^1 \rightarrow M$ of a 
non-reversible
Finsler metric on a differentiable manifold $M$ {\em geometrically
equivalent} if their traces $c_1(S^1)=c_2(S^1)$ coincide and 
if their orientations
coincide. The equivalence class is also called a {\em geometric
closed geodesic.}
For a closed geodesic $c_1$ there is a prime closed geodesic $c$
such that the set of all geometrically equivalent closed 
geodesics consists
of $u.c^m;m\ge 1, u\in S^1.$ 
Let $\theta: \Lambda M \rightarrow \Lambda M$ be the
orientation reversing, i.e. $(\theta c) (t)=c(1-t).$ 
For a non-reversible Finsler metric
this mapping in general does not leave the energy 
functional invariant.
And in general for a closed
geodesic $c$ the curve $\theta c$ is not a geodesic. 
\\
The second order behaviour of the energy functional in a neighborhood
of a closed geodesic is determined by its {\em index form} $\Hc $ which equals
the hessian $d^2E(c)$ 
of the energy functional by the second variational formula, cf. 
\cite[ch.2]{Ra4}. 
The {\em index} $\ind (c)$ of the closed geodesic $c$ is the index of the
index form $\Hc$ i.e. it is the maximal dimension of a subspace on which 
$\Hc$ is negative definite. The nullity $\nulli (c)$ is the nullity of the index
form $\Hc$ minus $1$. This convention is used since due to the $S^1$-action
the nullity of the index form $\Hc$ is at least $1.$ 
We call a Finsler metric {\em bumpy,}
if all closed geodesics are {\em non-degenerate,} 
i.e. for all closed geodesics $c$
the nullity $\nulli(c)=0$ vanishes.
\\
For a closed geodesic $c$ let
$\Lambda(c):=\{\gamma \in \Lambda M \,|\, E(\gamma)<E(c)\}.$ We call
\beqno
\overline{C}_*(c):=
H_*\left(\left(\Lambda(c)\cup S^1.c\right)/S^1,
\Lambda(c)/S^1;\q\right)
\eeqno
the {\em $S^1$-critical group} of $c,$ cf. \cite[ch.6.3]{Ra2}. 
Here $H_*$ is the singular homology, as coefficient field
we use in this paper the rationals $\q.$
Let $w_k(c):=\dim \overline{C}_k(c).$ 
\\
For the Morse theory of the energy functional it is 
important to describe the 
$S^1$-critical groups $\overline{C}_*(c^m)$ of the 
infinitely many critical points $S^1.c^m, m\ge 1$ 
produced by a single prime closed geodesic.
Therefore the sequence $\ind(c^m);m\ge1$ of the sequence 
of iterates $c^m$  is of great importance.
\begin{dfn}
\label{dfn:gam}{\rm (cf.~\cite[Def.1.6]{Ra1})}
For a closed geodesic $c$ 
let $\gamma_c \in \{\pm 1/2,\pm1\}$ be the 
invariant defined by
$\gamma_c >0$ if and only if $\ind (c)$ is even and
$|\gamma_c|=1$ if and only if $\ind (c^2)-\ind(c)$ is even. 
\end{dfn}
Then we obtain the following crucial 
\begin{lem}
\label{lem:local}{\rm (cf.~\cite[Prop.2.2]{Ra1})}
Let $c$ be a prime and non-degenerate closed geodesic, $m\ge 1.$
Then for all $k \ge 0:$
$$ w_k(c^m)= \left\{
\begin{array}{lcl}
1 &;& k=\ind(c^m), m \mbox{ \rm even or } \gamma_c=\pm1\\
&&\\
0 &;& \mbox{ \rm otherwise }
\end{array}
\right.
$$
\end{lem}
It follows from Bott's iteration formula \cite[thm. A,C]{Bo}
resp. Equation ~\ref{eq:bott} that the 
sequence $\ind(c^m);m\ge1$ grows almost linearly. 
Therefore the {\em average index}
\beqno
\alpha_c:=\lim_{m\to \infty}\frac{\ind(c^m)}{m}
\eeqno
is well-defined and one can show that
\begin{equation}
\label{eq:hb}
\left|\ind(c^m)-m\alpha_c\right|\le n-1
\end{equation}
for all $m\ge 1,$
where $n$ is the dimension of the manifold, cf. \cite[(1.4)]{Ra1}.
\\
Now we consider bumpy metrics with only a single
geometric closed geodesic.
Hence there is a prime closed geodesic
$c$ such that any closed geodesic on $M$ is of the form
$u.c^m$ for some $u \in S^1$ and $m\ge 1.$ 
Then we define for every
$k\ge0:$
\beq
\label{eq:wk}
w_k:=\sum_{m =1}^{\infty}w_k(c^m).
\eeq
This number gives the number of critical orbits $S^1.c^m$
whose $S^1-$critical group is non-trivial in dimension $k.$
The Morse inequalities relate these local invariants of the 
critical points to the Betti numbers
$$\betti_k:=b_k\left(\Lambda M/S^1,\Lambda^0M/S^1;\q\right)$$
of the pair $\left(\Lambda M/S^1;\Lambda^0 M/S^1)\right)$
of the quotient space $\Lambda M/S^1$ of the free loop space
$\Lambda M$  divided by
the $S^1$-action
and the space $\Lambda M^0:=\{c \in \Lambda M; E(c)=0\}$
of point curves. This space can be identified with the
manifold $M,$ it is the fixed point set of the $S^1-$action.
In particular we will use the following 
Betti numbers of the quotient space
$\Lambda S^n/S^1:$
\begin{pro}
\label{pro:betti3}{\rm (cf. \cite[Thm. 2.4]{Ra1})}
If $M$ is a simply-connected compact manifold rationally homotopy
equivalent to the $n$-dimensional sphere $S^n$ then the Betti
numbers $\betti_k:=b_k\left(\Lambda M/S^1,\Lambda^0M/S^1;\q\right)$
of the quotient $\Lambda M/S^1$ of the free loop space by
the $S^1-$action
satisfy:
\begin{enumerate}
\item $\betti_k\in \{0,1,2\},$ 
and $\betti_k\ge 1$ if and only if
$k \equiv n-1\pmod{2}$ and $k\ge n-1.$
\item If $n\equiv 0\pmod{2}$ then $\betti_k=2$ if and only if
$k=(2j+1)(n-1), j\ge 1.$
\item If $n\equiv 1\pmod{2}$ then $\betti_k=2$ if and only if
$k=j (n-1), j\ge 2.$
\end{enumerate}
\end{pro}
This Proposition follows from the form of the  Poincar\'e polynomials
$P(t)=\sum_k \betti_k t^k$ determined in \cite[(2.4),(2.5)]{Ra1}:
\begin{eqnarray*}
P(\Lambda S^n/S^1, \Lambda^0 S^n;\q)(t)=
\left\{
\begin{array}{ccc}
t^{n-1}\left\{\frac{1}{1-t^2}+
\frac{t^{2n-2}}{1-t^{2n-2}}\right\}
&;& n\equiv 0\pmod{2}\\
&&\\
t^{n-1}\left\{\frac{1}{1-t^2}+
\frac{t^{n-1}}{1-t^{n-1}}\right\}
&;& n\equiv 1\pmod{2}
\end{array}
\right.
\end{eqnarray*}
It follows that 
$B(n,1):=\lim_{N\to \infty}\sum_{j=0}^N (-1)^j b_j$ 
satisfies
\beq
\label{eq:bn}
B(n,1)=\left\{
\begin{array}{ccl}
- \frac{n}{2n-2}&;& n\equiv 0\pmod{2}\\
&&\\
\frac{n+1}{2n-2}&;& n\equiv 1\pmod{2}
\end{array}
\right.
\eeq
Then we collect the following conclusions from
the results in
\cite[(2.3),(2.6), Thm.3.1(a)]{Ra1} resp. \cite{Ra2}:
\begin{pro}
\label{pro:morse}
Let $F$ be a bumpy Finsler metric on a compact 
and simply-connected manifold of the rational homotopy type
of an $n-$sphere with
only  one geometric closed 
geodesic represented by the prime closed geodesic $c$ 
\begin{itemize}
\item[(a)]
The average index $\alpha_c$ and the invariant
$\gamma_c$ satisfy
\beq
\label{eq:al}
(-1)^{n-1} \,\frac{\alpha_c}{\gamma_c}=\,
 \left\{
\begin{array}{ccl}
2-\frac{2}{n}&;& n\equiv 0\pmod{2}\\
&&\\
2-\frac{4}{n+1}&;& n\equiv 1\pmod{2}
\end{array}
\right.
\,.
\eeq
\item[(b)] The numbers $w_k, k\ge 0$ 
of critical $S^1-$orbits whose $S^1-$critical group
is non-trivial in dimension $k$ (cf. Equation~\ref{eq:wk}) are 
bounded and there is a sequence $(q_k)_{k\ge 0}$ of non-negative
integers satisfying
\begin{equation}
\label{eq:morse-ineq}
w_k=\betti_k + q_k + q_{k-1}
\end{equation}
for all $k \ge 0.$
\end{itemize}
\end{pro}
\section{Bott's iteration formula and the index growth}
\label{sec:bott}
For a closed geodesic $c:\R \rightarrow M$ with
$c(t+1)=c(t)$ for all $t \in \R$ we define the {\em linearized Poincar\'e mapping}
$P_c:$ For $p=c(0)$ let $V$ be the 
$(n-1)-$dimensional orthogonal complement in the tangent space
$T_pM$ to the $1$-dimensional subspace generated by $c'(0)$ 
(with respect to the osculating Riemannian metric defined by the velocity
field $c'$.). Then let 
$$ P_c: V \oplus V \rightarrow V 
\oplus V\,;\, P_c(X,Y)=\left(J(1),\frac{\nabla}{dt} J(1)\right)$$
where $J$ is the Jacobi field determined uniquely by the
initial conditions $J(0)=X, \frac{\nabla}{dt}J(0)=Y$ along $c|[0,1].$
Here $\frac{\nabla}{dt}$ is the covariant derivative along $c.$ 
The canonical symplectic structure on $V\oplus V$ is preserved by $P_c,$
hence by choosing an orthonormal basis in $V$ we consider
$P_c$ as an element of the group $\symp (n-1)$ of linear symplectic
maps of $\R^{n-1}\oplus\R^{n-1},$ which is well-defined 
up to conjugation, i.e. independent of the
choice of $p=c(0)$ and the choice of an orthonormal basis of
$V.$
\\
We denote by 
$\spec(P) \subset \comp$ the set of eigenvalues of the 
complexification $P$
of $P_c$ and denote
by $V(z)=\ker \left(P - z\, Id\right)^{n-1}$
the generalized
eigenspace for $z \in \spec(P).$ Then 
\begin{equation}
\label{ncz1}
N_c(z):=\dim_{\comp}\left(P_c-z \,Id\right)
\le \dim_{\comp} V(z)\,,
\end{equation}
is called the {\em $z$-nullity} of $c.$
Since $P$ is symplectic, 
with $z \in \spec(P)$ also $z^{-1},\overline{z},
\overline{z}^{-1}\in \spec(P).$ 
If $c, c^2$ are non-degenerate (i.e. 
$1,-1 \not\in \spec(P)$) it follows:
\begin{equation}
\label{eq:ncz}
\sum_{z \in \spec(P), |z|=1, {\rm Im} z>0} N_c(z)\le n-1\,.
\end{equation}
In addition to the $z$-nullity $N_c(z)$
for a closed geodesic $c$ Bott defined the
$z$-index 
$I_c: \{z\in\comp;|z|=1\}  \rightarrow \n_0.$ 
Both functions have
the following properties, cf. \cite[thm. A,C]{Bo}, 
\cite[ch.4]{Ra2}:
The functions $I_c,N_c$ are invariant 
under conjugation,
i.e. $I_c(\overline{z})=I_c(z); N_c(\overline{z})=N_c(z)$ 
for all $z.$
The function $I_c$ is a continuous and constant function 
outside the
set $\spec(P_c)$ The finite jumps define the 
{\em splitting numbers}
$$
S_c^{\pm}(z)=\lim_{\phi \to \pm 0} 
I_c(z \exp\left(i \phi)\right)-I_c(z)
$$
satisfying 
\begin{equation}
\label{eq:sn}
0\le S_c^{\pm}(z)\le N_c(z)\,.
\end{equation}
The functions $I_c$ allow
the following formula due to Bott \cite[thm.A,C]{Bo} for the sequence $\ind(c^m)$ 
of the indices of the multiples $c^m$ of a closed geodesic $c:$
\begin{equation}
\label{eq:bott}
\ind(c^m)=\sum_{z^m=1} I_c(z)
\end{equation}
As an immediate consequence we obtain that for all
$m \ge1:$
\begin{equation}
\label{eq:indcm}
\ind(c^m) \ge \ind(c)\,.
\end{equation}
It was shown by Bott that the splitting numbers $S^{\pm}_c(z)$ 
only depend on the conjugacy class of the linearized Poincar\'e
mapping in the symplectic group. 
If $c$ is a closed geodesic of a bumpy metric and
$z_1=\exp(2\pi i t_1), \ldots, z_l=\exp (2\pi i t_l)$ 
are the
eigenvalues $z$ of the linearized Poinca\'e mapping with
$|z|=1, {\rm Im}(z)>0$ and $0=t_0<t_1<t_2<\ldots<t_l<t_{l+1}=1/2$ 
then $l \le n-1$ and the numbers $t_j$ are irrational
since $\nulli(c^m)=\sum_{z^m=1}N_c(z).$
With the help of the function $I_c$ we get the following
expression for the average index. Let 
$I_1=I_c(0)=I_c(\exp(2\pi i t)); t \in[0,t_1)$ and for
$j \in \{1,2,\ldots,l-1\}: 
I_j:=I_c\left(\exp\left(2\pi i t\right)\right); 
t \in \left(t_{j-1}, t_j \right)$
and $I_{l+1}:=I_c(-1)=I_c\left(\exp\left(2\pi i t\right)\right); 
t\in(t_l,1/2].$
Hence $\left\{I_c\left(\exp(z)\right)|z\in S^1-\spec (P_c)\right\}=
\left\{I_1,I_2,\ldots,I_{l+1}\right\}.$
Then Bott's formula Equation~\ref{eq:bott} implies
\begin{equation}
\label{eq:alI}
\alpha_c=\int_0^1 
I_c\left(\exp\left(2\pi i t\right)\right) \,dt=I_1 t_1
+ \sum_{j=1}^{l-1} I_j \left(t_j-t_{j-1}\right)
+I_l \left(\frac{1}{2}-t_l\right)
\end{equation}
Since
$ I_{j}-I_{j+1}=S^-(z_j)-S^+(z_j)$
we obtain from Equation~\ref{eq:sn}
$\left|I_j-I_{j+1}\right|\le N_c(z_j)$
where $N_c(z_j)\ge 1.$
Then Equation~\ref{eq:ncz} implies
\begin{equation}
\label{eq:nminus}
\sum_{j=1}^l \left|I_j-I_{j+1}\right|\le 
\sum_{j=1}^l N_c(z_j) \le n-1\,.
\end{equation}
The following Proposition will be crucial in the
Proof of our Theorem:
\begin{pro}
\label{pro:sn-1}
Let $c$ be a closed geodesic of a bumpy Finsler metric on
an $n$-dimensional manifold $M.$ We denote by
$z_j=\exp\left(2\pi i t_j\right), t_j\in (0,1/2),
j=1,2,\ldots,l\,;\, l\le n-1$ the 
eigenvalues
of the linearized Poincar\'e mapping $P_c$ 
whose imaginary part is positive with 
$0<t_1<t_2<\ldots<t_l<1/2.$ Let $t_0=0,t_{l+1}=1/2$
and $I_j:=I_c\left(\exp(2\pi i t)\right), 
t \in \left(t_{j-1},t_j\right); 
j=1,2,\ldots,l+1.$
\\
Let $\ind(c)=n-1, \ind\left(c^2\right)\ge n$ and 
$\alpha_c< 2 |\gamma_c|\,.$ \\
Then the invariants of the closed geodesic
$c$ satisfy:
\begin{itemize}
\item[(a)]
$\gamma_c=(-1)^{n-1};
\,\alpha_c >1\,;\,
\ind(c^2)=n+1$
\item[(b)]
$I_c(1)=I_1=n-1>I_2>\cdots>I_{l}=1;
I_{l+1}=I_c(-1)=2$ 
\item[(c)]
For all integers
$m\ge 1\,:\,\ind(c^{m+1}) \ge \ind(c^m)\,.$
\end{itemize}
\end{pro}
\begin{rem}
\label{rem:elliptic}
\rm
Since the number
$\overline{S}:=\sum_{j=1}^l|S^+(z_j)-S^-(z_j)|=
\sum_{j=1}^{n-1}|I_{j+1}-I_j|=n-1$ is 
maximal it follows from the formula for the 
splitting numbers depending on the symplectic 
normal form  given for example in
\cite[2.13]{BTZ} that the closed geodesic 
is of elliptic type,
i.e. the linearized Poincar\'e mapping has 
a symplectic normal form consisting only of
$2$-dimensional rotations.
\end{rem}
\begin{proof}
We define
\begin{equation*}
I_{min}:=\min \{I_1,I_2,\ldots,I_l\}\,;\,
r:=\min\left\{\left.j\in \{1,2,\ldots,l\}
\right|\,I_j=I_{\min}, 1\le j\le l\right\}
\end{equation*}
Since $\ind(c)=n-1$ we have $\gamma_c (-1)^{n-1}\in\{1/2,1\}.$
If $\gamma_c=(-1)^{n-1}/2$ then $\alpha_c <1$
by Proposition~\ref{pro:morse}(a),
hence $I_{\min}=0$ by Equation~\ref{eq:alI}. 
Then Equation~\ref{eq:nminus}
implies that $r=l; I_{l+1}=0.$ Bott's formula~\ref{eq:bott}
then implies that $\ind(c^2)=n-1+I_{l+1}=n-1$ contradicting the 
assumption $\ind(c^2) \ge n.$ \\
Hence we showed that $\gamma_c=(-1)^{n-1}.$
By Bott's formula~\ref{eq:bott} we conclude that
$\ind(c^2)=n-1+I_{l+1},$ hence by Definition~\ref{dfn:gam}
and the assumption $\ind(c^2) \ge n:$
$I_{l+1}$ is an even and positive number.
Since $\alpha_c<2$
Equation~\ref{eq:alI} implies that 
$I_{\min}=I_r \le 1.$ 
From $I_1=n-1$ we conclude
$I_r\in\{0,1\}.$
If $I_r=0$ then Equation~\ref{eq:nminus}
implies that $r=l+1$ and therefore 
$\ind(c^2)=n-1$ contradicting our assumption. \\
Hence we conclude $I_r=1.$ 
Equation~\ref{eq:nminus} implies that $r \in \{l,l+1\}$
and $I_{l+1}-I_l\le 1$ if $r=l.$
Since $I_{l+1}$ is even and positive we have $r=l$
and therefore $I_{l+1}=2.$
Since $I_1-I_{l}=n-2\,,\, I_{l+1}-I_{l}=1$ in
Equation~\ref{eq:nminus} the inequalities are
actually equalities and 
hence
$n-1=I_1>I_2>\cdots>I_{l}=1; I_{l+1}=2$
which finishes the proof of parts (a) and (b).
\\[1ex]
It remains to prove part (c):
We use the notation $e(x):=\exp(2\pi i x).$ 
Using Bott's formula~\ref{eq:bott} and
$I_c(-1)=2$ we obtain 
\beq
\label{eq:indcm+1}
\ind(c^{m+1}) - \ind(c^{m})=
A_m+B_m
\eeq
with
\beq
\label{eq:Am}
A_m:=\left\{
\begin{array}{ccc}
2 &;& m\equiv 1\pmod{2}\\
&&\\
2\,I_c\left(e\left(\frac{m}{2m+2}\right)\right)- \, 2 &;&
m\equiv 0\pmod{2}
\end{array}
\right.
\eeq
and
\beq
\label{eq:Bm}
B_m:= 2 \,\sum_{1\le j < m/2}\left\{
I_c\left(e\left(\frac{j}{m+1}\right)\right)
-I_c\left(e\left(\frac{j}{m}\right)\right)
\right\}
\eeq
From part (b) we conclude: 
$I_c\left(e\left(j/\left(m+1\right)\right)\right)
-I_c\left(e\left(j/m\right)\right)<0$ for some
$j \in \n, j<m/2$ implies
$j \in J_m:=\{p\in \n; p/(m+1)<t_l<p/m, p<m/2\}$
and therefore $I_c\left(e\left(j/\left(m+1\right)\right)\right)
-I_c\left(e\left(j/m\right)\right)=-1.$
On the other hand $J_m$ is either empty or 
consists of a single element $J_m=\{j_m\}.$
Hence $B_m \ge -2$ and $B_m=-2$ if and only if $J_m=\{j_m\}.$
Therefore we obtain immediately that 
$A_m+B_m\ge 0$ for
odd $m.$ Let $m$ be even, 
since $I_{min}=1$ we have $A_m\ge 0.$ On the other
hand $A_m=0$ implies that $m/(2m+2)<t_l<1/2.$
Hence
$J_m=\emptyset$ and $A_m+B_m\ge 0.$ If
$A_m >0$ then $A_m\ge 2,$ hence also in this case
$A_m+B_m \ge 0.$ 
Therefore we obtain in all cases 
$\ind\left(c^{m+1}\right)-\ind\left(c^m\right)=A_m+B_m\ge 0.$
\end{proof}
We use the {\em common index jump theorem} 
due to Y.Long and C.Zhu for a single closed geodesic
$c$ of a bumpy metric. For a closed geodesic $c$
of a bumpy metric the nullities
$\nulli(c^m)$ for all $m$ vanish as well as the
splitting number $S_c^+(1).$ Hence we obtain
as a particular case of the common index
jump theorem the following Proposition. Here we
note that the common index jump theorem also applies
to the closed geodesic problem on Finsler manifolds,
cf. \cite[Rem.12.2.5]{Lo}. 
\begin{pro}
\label{pro:jump}
{\rm \cite[Thm.4.3]{LZ}
\cite[Thm.11.2.1]{Lo}}
Let $c$ be a closed geodesic of a bumpy Finsler metric
with positive average index $\alpha_c.$ Then there are
infinitely many $k$ with
\begin{eqnarray*}
\ind\left(c^{2k+1}\right)-
\ind\left(c^{2k-1}\right)= 2 \,\ind(c)\,.
\end{eqnarray*}
\end{pro}
\section{Proof of the Theorem}
\label{sec:pr}
\begin{proof} {\em of the Theorem}
We assume that there is a prime closed geodesic $c$ such
that all closed geodesics on $M$ are geometrically equivalent
to $c$ and $n=\dim M \ge 3.$
For the average index $\alpha_c$ and the invariant
$\gamma_c$ 
we obtain from Proposition~\ref{pro:morse}(a):
\begin{equation}
\label{eq:alkl2}
1 \le \frac{\alpha_c}{|\gamma_c|} < 2\enspace; 
\enspace \alpha_c=|\gamma_c| \Leftrightarrow n=3\,.
\end{equation}
Let
$w_k:=\sum_{m=1}^{\infty} w_k(c^m),$
cf. Equation~\ref{eq:wk}. 
Since the Betti numbers
$\betti_k=b_k\left(
\Lambda M/S^1,\Lambda^0M/S^1;\q\right)$ vanish in dimensions
$k \equiv n \pmod{2},$ cf. Proposition~\ref{pro:betti3} 
we conclude from the
Morse inequalities Equation~\ref{eq:morse-ineq}
that the sequence $(q_j)_{j\ge0}$ vanishes identically, i.e.
\begin{equation}
\label{eq:bettiw}
w_k=\#\left\{m\in \n\,;\, \ind(c^m)=k, m \mbox{ odd or }
\gamma_c=\pm 1\right\} = \betti_k \in \{0,1,2\}
\end{equation}
for all $k \in \n_0$ 
with $\betti_k \in \{1,2\}$ if and only if $k \equiv n-1 \pmod{2}$
and $k \ge n-1\,,$ cf. Proposition~\ref{pro:betti3} 
and Lemma~\ref{lem:local}.
\\
We conclude from
Equation~\ref{eq:indcm} and from 
$\betti_0=\betti_1=\ldots=\betti_{n-2}=0, 
\betti_{n-1}=1$ (cf. Proposition~\ref{pro:betti3})
that 
\begin{equation}
\label{eq:indc}
\ind (c)=n-1\,.
\end{equation}
We conclude from Equation~\ref{eq:hb} and Equation~\ref{eq:alI}
that
$n-1 =\ind(c) \le \ind (c^2) < 2 \cdot 2 + n-1=n+3.$ 
If $\ind(c^2)=\ind(c)=n-1,$ then Definition~\ref{dfn:gam}
and Lemma~\ref{lem:local} imply
$w_{n-1}\ge w_{n-1}(c)+w_{n-1}(c^2)\ge 2$ contradicting
Equation~\ref{eq:bettiw}. Therefore $\ind(c^2)\ge n+1$
and we conclude from Proposition~\ref{pro:sn-1}(a):
$\gamma_c=(-1)^{n-1}$ and $\alpha_c>1.$ 
In dimension $n=3$ this yields a contradiction 
to Equation~\ref{eq:alkl2}
and finishes
the proof.
\\
Hence we can now assume $n\ge 4:$ 
Since the sequence $\ind(c^m),m\ge1$ is monotone increasing
by Proposition~\ref{pro:sn-1}(c) and 
$b_k=w_k\ge 1$ for all $k\equiv n-1\pmod{2}, k\ge n-1$ we 
have
$\ind(c^{m+1})-\ind(c^m)\in\{0,2\}$
for all $m\ge1.$ Therefore for $n\ge 4$ we obtain for
all $m\ge1:$
\beq
\label{eq:mplus}
\ind(c^{m+2})-\ind(c^m) \le 4\,.
\eeq
But we conclude from Proposition~\ref{pro:jump} that 
for $n\ge 4$ there is a $k_1\in \n$ with:
\beq
\ind(c^{2k_1+1})-\ind(c^{2k_1-1}) =2n-2\ge 6
\eeq
contradicting Equation~\ref{eq:mplus}. Therefore we obtain
a contradition to our assumption made at the beginning
of our proof
\end{proof}
\begin{rem}
\label{rem:pinching}
\rm
It is likely that two is not the optimal lower bound for the
number of closed geodesics of a non-reversible bumpy metric
on $S^n,n\ge 3.$ Katok's examples carry $2n$ closed geodesics
on $S^{2n}$ and $S^{2n-1}$. One can show that any non-reversible
and bumpy Finsler metric sufficiently $C^2-$close to the standard Riemannian
metric of constant curvature $1$
on $S^{2n}$ resp. $S^{2n-1}$ carries at least
$2n$ closed geodesics of length approximately $2\pi,$
cf. \cite[p.141]{Zi}. Estimates on the number of closed geodesics
under pinching assumptions for the flag curvature
depending on the reversibility are presented 
for example in \cite[Thm.3, Thm.8]{Ra2006}. Various questions resp.
conjectures about the optimal lower bound for an
arbitrary Finsler metric on $S^n$ are given
in \cite[p.155]{Zi} and \cite[ch.4]{Lo2}.
\end{rem}
\small

{\sc Universit\"at Leipzig, Mathematisches Institut\\ Augustusplatz 10/11,
D-04109 Leipzig, Germany}\\
{\tt rademacher@math.uni-leipzig.de}\\
{\tt www.math.uni-leipzig.de/\symbol{126}rademacher}
\end{document}